\documentclass[12pt,twoside]{article}    
\usepackage{amsfonts,amsmath,latexsym}

\newtheorem{theorem}{Theorem}[section]
\newtheorem{remark}[theorem]{Remark}
\newtheorem{lemma}[theorem]{Lemma}

\newtheorem{example}[theorem]{Example}
\numberwithin{equation}{section}
\newenvironment{proof}[1][Proof]{\textbf{#1.} }{\ \rule{0.5em}{0.5em}}

\makeatletter \@addtoreset{equation}{section} \makeatother
\begin{document}


\pagestyle{myheadings}

\markboth{\hfill {\small AbdulRahman Al-Hussein } \hfill}{\hfill
{\small Necessary conditions for optimality for SEEs} \hfill }


\thispagestyle{plain}


\begin{center}
{\large \textbf{Necessary conditions for optimality for stochastic evolution equations}$^{*}$\footnotetext{$^{*}$
This work is supported by the Science College Research Center at Qassim University, project no. SR-D-012-1610. }} \\
\vspace{0.7cm} {\large AbdulRahman Al-Hussein }
\\
\vspace{0.2cm} {\footnotesize
{\it Department of Mathematics, College of Science, Qassim University, \\
 P.O.Box 6644, Buraydah 51452, Saudi Arabia \\ {\emph E-mail:} alhusseinqu@hotmail.com}}
\end{center}

\begin{abstract}
This paper is concerned with providing the maximum principle for a control problem governed by a stochastic
evolution system on a separable Hilbert space. In particular, necessary conditions for
optimality for this stochastic optimal control problem are derived by using the adjoint backward
stochastic evolution equation. Moreover, all coefficients appearing in this system are allowed to
depend on the control variable. We achieve our results through the semigroup approach.
\end{abstract}

{\bf MSC 2010:} 60H10, 60H15, 93E20. \\

{\bf Keywords:} Stochastic evolution equation, optimal control, maximum principle, necessary conditions for optimality, backward stochastic evolution equation.

\sloppy

\section{Introduction}\label{sec1}
Consider a stochastic controlled problem governed by the following stochastic evolution equation (SEE):
\begin{eqnarray}\label{eq:first SEE}
 \left\{ \begin{array}{ll}
              d X(t) = ( A X(t) + b (X(t), \nu (t)) ) d t
              + \sigma (X(t), \nu (t)) d W(t) , \;\; t \in (0 , T] ,\\
             \; \, X(0) = x_0 .
         \end{array}
 \right.
\end{eqnarray}
We shall be interested in trying to minimize the cost functional, which is given by equation
(\ref{cost functional}) below, over a set of admissible controls.

This system is driven mainly by a possibly unbounded linear operator $A$ on a separable Hilbert space $H$ and a
cylindrical Wiener process $W$ on $ H .$ Here $\nu (\cdot )$ denotes a control process.

We shall derive the maximum principle for this control problem. More precisely, we shall concentrate on providing necessary conditions for
optimality for this optimal control problem, which gives this minimization.
For this purpose we shall apply the theory of \emph{backward stochastic evolution
equations} (BSEEs shortly) as in equation (\ref{adjoint-bse}) in Section~\ref{sec3}. These equations together with backward stochastic differential
equations (BSDEs) have become of great importance in a number of fields. For example in \cite{[Alh-COSA]}, \cite{[Alh-2012-1]}, \cite{[H-Pe96]}, \cite{[Oks05]}, \cite{[Oks-Zh05]}, \cite{[Pe-93]} and \cite{[Y-Z]} one can find applications of BSDEs to stochastic optimal control problems. Some of these references have also studied the maximum principle to find either necessary or sufficient conditions for optimality for stochastic differential equations (SDEs) or stochastic partial differential equations (SPDEs). Necessary conditions for optimality of the control process $\nu (\cdot )$ and its corresponding solution $X^{\nu (\cdot)}$ but for the case when the noise term $\sigma$ does not depend on $\nu (t)$ can be found in \cite{[H-Pe96]}.

In our work here we allow $\sigma$ to depend on the control variable and study a stochastic control problem associated with the former SEE. This control problem is explained in details in Section~\ref{sec2}, and the main theorem is stated in Section~\ref{sec3} and is proved together with all necessary estimates in Section~\ref{sec4}. Sufficient conditions for optimality for this optimal control problem can be found in \cite{[Alh-2012-1]}. We refer the reader also to \cite{[Alh-COSA]}.

On the other hand, we recall that control problems governed by SPDEs that are driven by martingales are studied in \cite{[Alh-AMO-10]}. In fact in \cite{[Alh-AMO-10]} we derived the maximum principle (necessary conditions) for optimality of stochastic systems governed by SPDEs. The technique used there relies heavily on the variational approach. The reason beyond that is that the only known way until now to find solutions to the resulting adjoint BSPDEs is achieved through the same variational approach, and is established in details in \cite{[Alh-Stoc09]}. Thus the semigroup approach to get mild solutions (as done here in Theorem~\ref{th:solution of adjointeqn} below and in Section~\ref{sec3}) cannot be used to study such adjoint BSPDEs considered in \cite{[Alh-AMO-10]}. Moreover, it is not obvious how one can allow the control variable $\nu(t)$ to enter in the noise term and in particular in the mapping $G $ in equation (1.1) of \cite{[Alh-AMO-10]} and obtain a result like Theorem~\ref{main thm} below. This problem is still open and is also pointed out in \cite[Remark~6.4]{[Alh-AMO-10]}.

In the present work, we shall show how to handle this open problem in great success, and as we stated earlier, we can and will allow all coefficients in (\ref{eq:first SEE}) and especially in the diffusion term to depend on the control variable $\nu (t) .$ We emphasize that our work here does not need go through the technique of Hamilton-Jacobi-Bellman equations nor the technique of viscosity solutions. We refer the reader to \cite{[Deb-H-Tess011]} for this business and to \cite{[Cerrai_book]} and some of the related references therein for the semi-group technique. Thus our results here are new. In this respect we thank the anonymous referee for pointing out the recent and relevant work of Fuhrman et al. in \cite{Fuh-Tes012}.

\section{Statement of the problem}\label{sec2}
Let $(\Omega, \mathcal{F}, \mathbb{P})$ be a complete probability space and denote
by $\mathcal{N}$ the collection of $\mathbb{P}$\,-\,null sets of
$\mathcal{F} .$ Let $\{ W (t) ,\, 0 \leq t \leq T \} $ be a
cylindrical Wiener process on $H$ with its completed natural
filtration ${\mathcal{F}_t = \sigma \{ \ell \circ W (s) \, , \; 0
\leq s \leq t \, , \ell \in H^{*} \}\vee \mathcal{N} } , \; t \geq 0
;$ see \cite{[preprint_MRT]} for more details.

For a separable Hilbert space $E$ denote by $L^2_{\mathcal{F}} (
0, T; E )$ to the space of all $\{ \mathcal{F}_t , 0
\leq t \leq T \}$\,-\,progressively measurable processes $f$
with values in $E$ such that \[ \mathbb{E}\; [ \int_0^T |
f(t) |_{E}^2 \; dt ] < \infty .\] This space is Hilbert with respect to the
norm \[ || f || \; = \Big{(} \mathbb{E}\; [ \int_0^T |
f(t) |_{E}^2 \; dt ] \Big{)}^{1/2}\, . \]
Moreover, if $ f \in L^2_{\mathcal{F}} (
0, T; L_2 (H) ) ,$ where $L_2 (H)$ is the space of all
Hilbert-Schmidt operators on $H ,$  the stochastic integral $\int f(t)
dW(t)$ can be defined and is a continuous stochastic martingale in $H .$ The norm and inner product on $L_2 (H)$ will be denoted respectively by $|| \cdot ||_2$ and  $\big{<} \cdot , \cdot \big{>}_2 .$

\bigskip

Let us assume that $\mathcal{O}$ is a separable Hilbert space equipped with an inner
product $\big{<} \cdot , \cdot \big{>}_{\mathcal{O}}$, and $U$ is a convex subset of $\mathcal{O} .$ We say that
$\nu (\cdot ) : [0 , T]\times \Omega \rightarrow \mathcal{O}$ is \emph{admissible}
if $\nu (\cdot ) \in  L^2_{\mathcal{F}} ( 0 , T ; \mathcal{O} )$ and $\nu (t) \in U \; \; a.e., \; a.s.$
The set of admissible controls will be denoted by $\mathcal{U}_{ad} .$

Suppose that $b: H \times \mathcal{O} \rightarrow H$ and $\sigma: H \times \mathcal{O} \rightarrow L_2(H)$ are two continuous mappings, and consider the following controlled SEE:
\begin{eqnarray}\label{forward-see}
\left\{ \begin{array}{ll}
              d X(t) = ( A X(t) + b (X(t), \nu (t)) ) d t
              + \sigma (X(t), \nu (t))  d W(t) ,\\
             \; \, X(0) = x_0 ,
         \end{array}
 \right.
\end{eqnarray}
where $\nu (\cdot ) \in \mathcal{U}_{ad} .$
A solution (in the sense of the following theorem) of (\ref{forward-see}) will be denoted by $X^{\nu (\cdot )}$ to indicate
the presence of the control process $\nu (\cdot) .$

Let $\ell : H \times \mathcal{O} \rightarrow \mathbb{R}$ and $\phi : H \rightarrow \mathbb{R}$ be two measurable mappings such that the following {\it cost functional} is defined:
\begin{equation}\label{cost functional}
J(\nu (\cdot ) ) : = \mathbb{E} \; [ \; \int_0^T \ell ( X^{\nu (\cdot )}
(t) , \nu (t) )  dt + \phi ( X^{\nu (\cdot )} (T) ) \; ] , \;\; \nu (\cdot ) \in \mathcal{U}_{ad} .
\end{equation}
For example one can take $\ell$ and $\phi$ to satisfy the assumptions of Theorem~\ref{main thm} in Section~\ref{sec3}.

\bigskip

The optimal control problem of the system (\ref{forward-see}) is to find the {\it
value function } $$J^{*} : = \inf \{ J(\nu (\cdot ) ) : \; \nu (\cdot ) \in \mathcal{U}_{ad} \}$$ and an \emph{optimal control} $\nu^{*} (\cdot
) \in \mathcal{U}_{ad} $ such that
\begin{equation}\label{control problem}
J^{*} = J(\nu^{*} (\cdot ) )  .
\end{equation}
If this happens, the corresponding solution $X^{\nu^{*} (\cdot )}$ is
called \emph{an optimal solution} of the stochastic control problem
(\ref{forward-see})--(\ref{control problem}) and
$( X^{\nu^{*} (\cdot )}\, , \nu^{*} (\cdot ) )$ is called an \emph{optimal pair.}

We close this section by the following theorem.
\begin{theorem}\label{thm1}
Assume that $A$ is an unbounded linear operator on $H$ that
generates a $C_0$-semigroup $\{ S(t) ,\; t \geq 0 \} $ on $H$, and
$b , \sigma$ are continuously Fr\'echet differentiable with respect to
$x$ and their derivatives $b_x \, , \,  \sigma_x$ are uniformly bounded.
Then for every $\nu (\cdot ) \in \mathcal{U}_{ad} $ there exists a unique mild solution $X^{\nu (\cdot )}$
on $[0, T] $ to (\ref{forward-see}). That is $X^{\nu (\cdot )}$ is a progressively measurable stochastic
process such that $X (0) = x_0 $ and for all $t \in [0 , T ] ,$
\begin{eqnarray}\label{thm1:solution}
X^{\nu (\cdot )} (t) &=& S ( t ) x_0 + \int_0^t S ( t - s) b(X^{\nu (\cdot )} (s) , \nu (s) ) ds \nonumber  \\
& & \hspace{2cm} + \, \int_0^t S ( t - s) \, \sigma (X^{\nu (\cdot )} (s) , \nu (s) ) \, d W(s) .
\end{eqnarray}
\end{theorem}
The proof of this theorem can be derived in a similar way to those
in \cite[Chapter 7]{[Da-Za]} or \cite{[Ichi]}.

\smallskip

From here on we shall assume that $A$ is the infinitesimal generator of a $C_0$-semigroup $\{ S(t) ,\; t \geq 0 \} $ on $H .$
Its adjoint operator $A^{*} : \mathcal{D} ( A^{*} ) \subset H \rightarrow H$ is then the infinitesimal generator of
the adjoint semigroup $\{ S^{*} (t) ,\; t \geq 0 \} $ of $\{ S (t) \, , t \geq 0 \} . $

\smallskip

\section{Stochastic maximum principle}\label{sec3}
It is known from the literature that BSDEs play
a fundamental role in deriving the maximum principle for SDEs. In this section we shall search for such a role for SEEs like
(\ref{forward-see}). To prepare for this business let us first define the \emph{Hamiltonian} by the following formula:
$$ \mathcal{H} : H\times \mathcal{O} \times H
\times L_2(H) \rightarrow \mathbb{R} ,$$
\begin{eqnarray}\label{def:Hamiltonian}
\mathcal{H} (x, \nu , y , z ): = \ell (x , \nu) \, + \big{<} b(x , \nu ) , y\big{>}_{H} +
\big{<}\sigma (x , \nu ) , z \big{>}_2 .
\end{eqnarray}
Then we consider the following BSEE on $H$:
\begin{eqnarray}\label{adjoint-bse}
\left\{ \begin{array}{ll}
             -\, d Y^{\nu (\cdot )} (t) = & \big(\, A^{*} \, Y^{\nu (\cdot )} (t) + \nabla_{x} \mathcal{H}
             ( X^{\nu (\cdot )}(t), \nu (t), Y^{\nu (\cdot )}(t) , Z^{\nu (\cdot )} (t) ) \,
             \big)\, dt \\& \hspace{1.70in} - Z^{\nu (\cdot )} (t) d W(t) , \; \; 0 \leq t < T,  \\
             \; \; \; Y^{\nu (\cdot )} (T) = & \nabla \phi (X^{\nu (\cdot )}(T)),
         \end{array}
 \right.
\end{eqnarray}
where $\nabla \phi$ denotes the gradient of $\phi ,$ which is
defined, by using the directional derivative $D\phi (x) (h)$ of
$\phi$ at a point $x \in H$ in the direction of $h \in H ,$ as
$\big{<} \nabla \phi (x) , h \big{>}_{H} = D\phi (x) (h) \; ( \, = \phi_x (h) \, ).$
This equation is the adjoint equation of (\ref{forward-see}).
\bigskip

As in the previous section a \emph{mild solution} (or a solution) of (\ref{adjoint-bse}) is a pair $( Y , Z ) \in L^2_{\mathcal{F}}( 0,
T; H ) \times L^2_{\mathcal{F}}( 0, T; L_2(H)) $ such that we have
$\mathbb{P}$\,-\,a.s. for all $t \in [0 , T]$
\begin{eqnarray}\label{bse:sol}
Y^{\nu (\cdot )} (t) &=& S^{*} ( T - t )\,  \nabla \phi (X^{\nu (\cdot )}(T)) \nonumber \\ &
& + \, \int_t^T S^{*} ( s - t )\, \nabla_{x} \mathcal{H} (
X^{\nu (\cdot )}(s), \nu (s), Y^{\nu (\cdot )}(s) , Z^{\nu (\cdot )} (s) ) ds \nonumber \\
& & \hspace{1.5in} - \, \int_t^T S^{*} ( s - t)\, Z^{\nu (\cdot )}(s) dW(s) .
\end{eqnarray}
\begin{theorem}\label{th:solution of adjointeqn}
Assume that $b , \sigma , \ell , \phi$ are continuously Fr\'echet differentiable with respect to
$x ,$ the derivatives $b_x , \sigma_x , \sigma_{\nu} , \ell_x $ are uniformly bounded, and
$$ | \phi_x |_{L(H,H)} \leq k \, ( 1 + |x|_{H} ) $$ for some constant $k > 0 .$

Then there exists a unique (mild) solution $( Y^{\nu (\cdot )} , Z^{\nu (\cdot )})$ of BSEE~(\ref{adjoint-bse}).
\end{theorem}

The proof of this theorem can be found in \cite{[BSEEs]} or \cite{[H-Pe91]}. An alternative proof by using finite dimensional framework through the Yosida approximation of $A$ can be found in \cite{[Tess96]}.

\bigskip

Our main result is the following.
\begin{theorem}\label{main thm}
Suppose that the following two conditions hold. \\
(i)\; $b , \sigma , \ell$ are continuously Fr\'echet differentiable with respect to
$x , \nu ,$ $\phi$ is continuously Fr\'echet differentiable with respect to
$x ,$ the derivatives $b_x , \, b_{\nu} , \, \sigma_x , \sigma_{\nu} , \ell_x , \, \ell_{\nu}$ are uniformly bounded, and
$$ | \phi_x |_{L(H,H)} \leq k \, ( 1 + |x|_{H} ) $$ for some constant $k > 0 .$
\\
(ii)\; $\ell_x$ is Lipschitz with respect to $u $ uniformly in $x .$

If $( X^{\nu^{*} (\cdot )} , \nu^{*} (\cdot ) )$ is an optimal
pair for the control problem (\ref{forward-see})--(\ref{control problem}), then there exists a unique solution $( Y^{\nu^{*} (\cdot )},
Z^{\nu^{*} (\cdot )} )$ to the corresponding BSEE~(\ref{adjoint-bse}) s.t. the following inequality holds:
\begin{eqnarray*}\label{ineq1:main}
&& \big{<} \, \nabla_{\nu } \mathcal{H} (X^{\nu^{*} (\cdot )}(t), \nu^{*} (t),
Y^{\nu^{*} (\cdot )}(t) , Z^{\nu^{*} (\cdot )}(t) ) \, , \, \nu^{*} (t) - \nu \, \big{>}_{\mathcal{O}} \leq 0  \\ \smallskip
&& \hspace{3.25in}  \text{a.e.} \; t \in [0 , T],\; \text{a.s.} \; \forall \; \nu \in U .  \nonumber
\end{eqnarray*}
\end{theorem}

\bigskip

The proof of this theorem will be given in Section~\ref{sec4} below. Now to illustrate this theorem let us present an example.
\begin{example}\label{ex1}
Let $H$ and $\mathcal{O}$ be two separable Hilbert spaces as considered earlier, and let $U = \mathcal{O} .$ We shall study in this example a special case of the control problem (\ref{forward-see})--(\ref{control problem}). In particular, given $\phi$ as in Theorem~\ref{main thm}, we would like to minimize the cost functional:
\begin{equation}\label{ex1:eq1}
J( \nu (\cdot ) ) = \mathbb{E} \; [ \; \int_0^T  | \nu (t)
|^2_{\mathcal{O}} \; dt \; ] + \mathbb{E} \; [ \; \phi ( X^{\nu (\cdot)}  (T) ) \; ]
\end{equation}
subject to:
\begin{eqnarray}\label{ex1:f-see}
 \left\{ \begin{array}{ll}
              d X^{\nu (\cdot)} (t) = ( \, A \, X^{\nu (\cdot)} (t) + B \, \nu (t) \, )\, d t
              + \, D \, \nu (t) \, d W(t) , \;\; t \in (0 , T] ,\\
             \; X^{\nu (\cdot)} (0) = x_0 \in H ,
         \end{array}
 \right.
\end{eqnarray}
where $B$ is a bounded linear operator from $\mathcal{O}$ into $H$ and $D$ is another bounded linear operator from $\mathcal{O}$ into $L_2(H) .$

The Hamiltonian is then given by the formula:
\begin{equation*}
 \mathcal{H} ( x , \nu , y , z ) = | \nu |^2_{\mathcal{O}} + \big<
B \, \nu  \, , y \big>_{H} + \big< D \nu \, , z \big>_{L_2 (H)} \, ,
\end{equation*}
where $( x , \nu , y , z ) \in H \times \mathcal{O} \times H \times L_2 (H) ,$ and the adjoint BSEE is
\begin{eqnarray}\label{eq:ex1-bsee}
 \left\{ \begin{array}{ll}
             -\, d Y^{\nu (\cdot)}  (t) = & A^{*} \, Y^{\nu (\cdot)} (t) dt - Z^{\nu (\cdot)} (t) d W(t) , \; \; \; t \in [0 , T) , \\
             \hspace{0.45cm}  Y^{\nu (\cdot)}  (T) = & \nabla \phi (X^{\nu (\cdot)} (T)) .
         \end{array}
 \right.
\end{eqnarray}

From the construction of the solution of (\ref{eq:ex1-bsee}), as e.g. in \cite[Lemma~3.1]{[BSEEs]}, this BSEE attains an explicit solution:
\[ Y^{\nu (\cdot)}  (t) = \mathbb{E} \, [ \; S^* (T-t) \, \nabla \phi (X^{\nu (\cdot)} (T)) \; | \; \mathcal{F}_t \; ] , \]
\[ Z^{\nu (\cdot)}  (t) = S^* (T-t) \, R^{\nu (\cdot)}  (t) , \]
where $R^{\nu (\cdot)} $ is the unique element of $L^2_{\mathcal{F}} (
0, T; L_2 (H) )$ satisfying
\[ \nabla \phi (X^{\nu (\cdot)} (T)) = \mathbb{E} \, [ \; \nabla \phi (X^{\nu (\cdot)} (T)) \; ] + \int_0^T R^{\nu (\cdot)}  (t) \, dW(t) . \]

On the other hand, for fixed $(x , y , z) ,$ we note that the function ${ \nu \mapsto \mathcal{H} ( x , \nu , y , z ) }$ attains its minimum at
$\nu = \frac{1}{2} \, \big( B^{*} \, y + D^* \, z \big) \; ( \, \in U \, ) ,$ where $B^*: H \rightarrow \mathcal{O}$ and $D^*: L_2(H) \rightarrow \mathcal{O}$ are the adjoint operators of $B$ and $D$ respectively. So we elect
\begin{equation}\label{eq:candidate optimal control}
\nu^* (t , \omega ) = \, \frac{1}{2} \, \big( B^{*} \, Y^{\nu^*(\cdot)}(t , \omega ) + D^* \, Z^{\nu^*(\cdot)}(t , \omega ) \big)
\end{equation}
as a candidate optimal control.

It is easy to see that with these choices all the requirements of Theorem~\ref{main thm} are verified. Hence this candidate $\nu^* (\cdot )$ given in (\ref{eq:candidate optimal control}) is an optimal control for the problem (\ref{ex1:eq1})--(\ref{ex1:f-see}), and its corresponding optimal solution $X^{\nu^{*} (\cdot )}$ is the solution of the following SEE:
\begin{eqnarray*}\label{ex1:optimal-solution}
 \left\{ \begin{array}{ll}
              d X^{\nu^*(\cdot)}(t) & = \Big( A \; X^{\nu^*(\cdot)}(t) + \, \frac{1}{2} \, B \, \big( \, B^{*} \, Y^{\nu^*(\cdot)}(t) + D^* \, Z^{\nu^*(\cdot)}(t) \, \big) \Big) d t \\
             &  \hspace{2.25cm} + \, \frac{1}{2} \,  D \, \big( \, B^{*} \, Y^{\nu^*(\cdot)}(t) + D^* \, Z^{\nu^*(\cdot)}(t) \, \big) \Big) d W(t) , \;\; t \in (0 , T] ,\\
             \; X^{\nu^*(\cdot)} (0) &= x_0 .
         \end{array}
 \right.
\end{eqnarray*}

Finally, the value function attains the formula
\begin{eqnarray*}
J^{*} &=& \frac{1}{4} \; \mathbb{E} \, \big[ \, \int_0^T  | \, B^{*} \, Y^{\nu^*(\cdot)}(t) + D^* \, Z^{\nu^*(\cdot)}(t) |^2_{\mathcal{O}} \, dt \, \big] + \mathbb{E} \, [ \, \phi (X^{\nu^*(\cdot)}(T)) \, ] .
\end{eqnarray*}
\end{example}

\smallskip

\begin{remark}
A concrete example in the setting of Example~\ref{ex1} can be constructed by taking $H = \mathcal{O} = L^2 ( \mathbb{R}^d ) , \; d \geq 1 ,$
$A = \frac{1}{2} \Delta$ (half-Laplacian), $B = id_{H} ,$ $D \nu := \big< v \, , h \big>_{H} \, \mathcal{Q}^{1/2} , $ $\phi (x) = \big{<} \rho \, , x \big{>}_{H} ,$ for some fixed elements $h , \rho$ of $H $ and a positive definite nuclear operator $\mathcal{Q}$ on $H .$

The computations in this case of $\mathcal{H} , Y^* , Z^* , \nu^* , X^*$ become direct from the corresponding equations in Example~\ref{ex1}.
\end{remark}

\section{Proofs}\label{sec4}
Let $\nu^{*} (\cdot )$ be an optimal control and $X^{*} \equiv X^{\nu^{*} (\cdot ) }$ be the corresponding solution of (\ref{forward-see}). Let $\nu (\cdot) $ be an element of $L^2_{\mathcal{F}}( 0, T; \mathcal{O})$ such that $\nu^{*} (\cdot ) + \nu (\cdot ) \in \mathcal{U}_{ad} .$
For a given $0 \leq \varepsilon \leq 1$ consider the variational control: $$\nu_{\varepsilon } (t) =  \nu^{*} (t) + \varepsilon \, \nu(t) , \; \; t \in [0 , T] .$$
We note that the convexity of $U$ implies that $\nu_{\varepsilon } (\cdot ) \in \mathcal{U}_{ad} .$
Considering this control $ \nu_{\varepsilon } (\cdot )$ we shall let $X^{\nu_{\varepsilon } (\cdot ) }$ be the solution of the SEE~(\ref{forward-see}) corresponding to $ \nu_{\varepsilon } (\cdot ),$ and denote it briefly by $X_{\varepsilon } .$

Let $p$ be the solution of the following linear equation:
\begin{eqnarray}\label{eq:p}
 \left\{ \begin{array}{ll}
               d p(t) = \big( A \, p(t) + b_x (X^{*}(t) , \nu^{*} (t) ) \, p (t)
             + b_{\nu } (X^{*}(t) , \nu^{*} (t) ) \, \nu (t) \big) dt \\ \hspace{1in} + \, \big( \, \sigma_x (X^{*}(t) , \nu^{*} (t) ) \,
               p (t) + \sigma_{\nu } (X^{*}(t) , \nu^{*} (t) ) \, \nu (t) \, \big) d W(t) ,\\
               \; p (0) = 0 .
         \end{array}
 \right.
\end{eqnarray}

\bigskip

The following three lemmas contain estimates that will play a vital role in deriving the desired variational equation and the maximum principle for our control problem.
\begin{lemma}\label{lem:1st estimate}
Assume condition (i) of Theorem~\ref{main thm}. Then
\[
\sup_{t \in [0 , T]} \mathbb{E} \, [ \, | p (t) |^2 \, ] \, < \infty .
\]
\end{lemma}
\begin{proof}
The solution of (\ref{eq:p}) is given by the formula
\begin{eqnarray}\label{eq:p-mild formula}
&& \hspace{-1.25cm} p(t) =\int_0^t S ( t - s) \big( b_x (X^{*}(s) , \nu^{*} (s) ) \, p (s) + b_{\nu } (X^{*}(s) , \nu^{*} (s) ) \, \nu (s) \big) ds \nonumber  \\
& & \hspace{-0.65cm} + \, \int_0^t S ( t - s) \,  \big( \, \sigma_x (X^{*}(s) , \nu^{*} (s) ) \,
               p (s) + \sigma_{\nu } (X^{*}(s) , \nu^{*} (s) ) \, \nu (s) \, \big) \, d W(s) .
\end{eqnarray}
By using Minkowski's inequality (triangle inequality), Holder's inequality, Burkholder's inequality for stochastic convolution together with assumption (i) and Gronwall's inequality we obtain easily
\begin{equation}\label{eq:est-p}
\sup_{t \in [0 , T]} \, \mathbb{E}\, [ \, |\, p (t) \, |^2 \, ] \leq \, C
\end{equation}
for some constant $C > 0 .$
\end{proof}

\begin{lemma}\label{lem:2nd estimate}
Assuming condition (i) of Theorem~\ref{main thm}, we have
\[
\sup_{t \in [0 , T]} \mathbb{E} \, [ \, | X_{\varepsilon} (t) - X^{*} (t) |^2 \, ] \, = \, O (\varepsilon^2 ) .
\]
\end{lemma}
\begin{proof}
Observe first from (\ref{thm1:solution}) that
\begin{eqnarray}\label{eq1::difference of x}
&& \hspace{-1.25cm} X_{\varepsilon} (t) - X^{*} (t) = \int_0^t S ( t - s) \big( \, b(X_{\varepsilon} , \nu_{\varepsilon} (s) ) - b(X^{*} (s) , \nu^{*} (s) ) \, \big) ds \nonumber  \\
& & \hspace{1.25cm} + \, \int_0^t S ( t - s) \, \big( \, \sigma (X_{\varepsilon} , \nu_{\varepsilon} (s) ) - \sigma (X^{*} (s) , \nu^{*} (s) ) \, \big) d W(s) .
\end{eqnarray}
Hence
\begin{eqnarray}\label{ineq2::difference of x}
&& \mathbb{E} \, [ \, | X_{\varepsilon} (t) - X^{*} (t) |^2 \, ] \leq 2 M^2 \, T \; \mathbb{E}\,  [ \, \int_0^t | \, b( X_{\varepsilon }(s) , \nu_{\varepsilon } (s) ) - b( X^{*}(s) , \nu^{*} (s) ) \, |^2  \; ds \, ] \nonumber \\
&& \hspace{1.25in} +\; 2 M^2 \; \mathbb{E}\,  [ \, \int_0^t || \sigma (X_{\varepsilon} , \nu_{\varepsilon} (s) ) - \sigma (X^{*} (s) , \nu^{*} (s) ) ||^2_{2}  \; ds \, ] ,
\end{eqnarray}
where $M : = \displaystyle{\sup_{t\in [0 , T]}} || S(t) ||_{L(H,H)} .$

Secondly, from condition (i) we get
\begin{eqnarray}\label{ineq3::difference of x}
&& \hspace{-1cm} \mathbb{E}\,  [ \, \int_0^t | \, b( X_{\varepsilon }(s) , \nu_{\varepsilon } (s) ) - b( X^{*}(s) , \nu^{*} (s) ) \, |^2  \; ds \, ] \nonumber \\ &\leq& 2 \, \mathbb{E}\,  [ \, \int_0^t | \, b( X_{\varepsilon }(s) , \nu_{\varepsilon } (s) ) - b( X^{*}(s) , \nu_{\varepsilon } (s) ) \, |^2  \; ds \, ]
\nonumber \\ && + \; 2 \,  \mathbb{E}\,  [ \, \int_0^t | \, b( X^{*}(s) , \nu_{\varepsilon } (s) ) - b( X^{*}(s) , \nu^{*} (s) ) \, |^2  \; ds \, ]
\nonumber \\ &=& 2 \, \mathbb{E}\,  [ \, \int_0^t  |\, \tilde{b}_x (s , \varepsilon ) ( X_{\varepsilon} (s) - X^{*} (s) ) |^2 \, ds \, ] + \, 2 \, \mathbb{E}\,  [ \, \int_0^t | \delta_{\varepsilon } b (s )   |^2 \; ds \, ]
\nonumber \\ &\leq&  2\, C_1 \, \mathbb{E}\,  [ \, \int_0^t | \, X_{\varepsilon} (s) - X^{*} (s) \, |^2 \, ds \, ] + 2 \,  C_2 \; \varepsilon^2 ,
\end{eqnarray}
where, for $y \in H ,$
\[ \tilde{b}_x (s , \varepsilon ) (y) = \int_0^1 b_x (  X^{*}(s) + \theta ( X_{\varepsilon} (s) - X^{*} (s) ) , \nu_{\varepsilon } (s) ) (y) d\theta  , \]
\[ \delta_{\varepsilon } b (s) = b (X^{*}(s) , \nu_{\varepsilon } (s) ) - b (X^{*}(s) , \nu^{*} (s) ) , \]
$C_1$ is a positive constant, and $C_2$ is another positive constant coming thanks to (i) from the following inequality:
\begin{eqnarray}\label{ineq4::difference of x}
&& \hspace{-0.7cm} \mathbb{E}\,  [ \, \int_0^T | \, \delta_{\varepsilon } b (s) \, |^2 \, ds \, ]
=
\mathbb{E}\,  [ \, \int_0^T | \, b (X^{*}(s) , \nu_{\varepsilon } (s) ) - b (X^{*}(s) , \nu^{*} (s) ) \, |^2 \, ds \, ]
\nonumber \\
&& \hspace{-0.5cm} =
\mathbb{E}\,  [ \, \int_0^T | \, \int_0^1 b_{\nu } (  X^{*}(s) , \nu^{*} (s) + \theta ( \nu_{\varepsilon} (s) - \nu^{*} (s) ) ) \, ( \nu_{\varepsilon} (s) - \nu^{*} (s) ) \, d\theta \, |^2 \, ds \, ]
\nonumber \\
&& \hspace{4in}
\leq  C_2 \, \varepsilon^2 .
\end{eqnarray}

Similarly,
\begin{eqnarray}\label{ineq5::difference of x}
&& \mathbb{E}\,  [ \, \int_0^t | \, \sigma( X_{\varepsilon }(s) , \nu_{\varepsilon } (s) ) - \sigma( X^{*}(s) , \nu^{*} (s) ) \, |^2  \; ds \, ] 
\nonumber \\ &&  \hspace{1in} \leq  2\, C_3 \, \mathbb{E}\,  [ \, \int_0^t | \, X_{\varepsilon} (s) - X^{*} (s) \, |^2 \, ds \, ] + 2 \,  C_4 \; \varepsilon^2 ,
\end{eqnarray}
for some positive constants $C_3, \, C_4 .$

Finally, by applying (\ref{ineq3::difference of x}), (\ref{ineq5::difference of x}) in (\ref{eq1::difference of x}) and then using Gronwall's inequality we find that
\begin{eqnarray}\label{ineq5:lem2}
& & \hspace{-1cm} \mathbb{E}\, [ \, |\, X_{\varepsilon} (t) - X^{*} (t) \, |^2 \, ] \leq C_5 \, \varepsilon^2
\end{eqnarray}
for some constant $C_5 > 0$ that depends in particular on $C_i , \, i = 1 , \ldots , 4,$ and $M .$ Hence the proof is complete.
\end{proof}

\bigskip

Keeping the notations $\tilde{b}_x$ and $\delta_{\varepsilon} b$ used in the preceding proof let us state the following lemma.
\begin{lemma}\label{lem:3rd estimate}
Let $\eta_{\varepsilon} (t) = \frac{X_{\varepsilon} (t) - X^{*} (t)}{\varepsilon } - p (t) .$ Then, under condition (i) of Theorem~\ref{main thm},
\[  \lim_{\varepsilon \rightarrow 0^{+} } \; \sup_{t \in [0 , T]} \mathbb{E} \, [ \, | \eta_{\varepsilon} (t) |^2 \, ] \, = \, 0 . \]
\end{lemma}
\begin{proof}
From the corresponding equations (\ref{forward-see}) and (\ref{eq:p}) we deduce that
\begin{eqnarray}\label{eq:eta}
&& \hspace{-0.75cm} \eta_{\varepsilon} (t) =  \nonumber \\ &&  \int_0^t S(t-s) \big[ \, \frac{1}{\varepsilon }\; \big( b( X_{\varepsilon}(s) , \nu_{\varepsilon } (s) ) - b( X^{*}(s) , \nu_{\varepsilon} (s) ) \big) -  b_x ( X^{*}(s) , \nu^{*} (s) ) \, p (s) \, \big] \, ds \nonumber \\
& & + \int_0^t S(t-s) \big[ \, \frac{1}{\varepsilon }\; \delta_{\varepsilon } b (s) - b_{\nu } ( X^{*}(s) , \nu^{*} (s) ) \, \nu (s) \, \big] \, ds  \nonumber \\
& & +  \int_0^t S(t-s) \big[ \, \frac{1}{\varepsilon }\; \big( \sigma (X_{\varepsilon }(s) , \nu_{\varepsilon } (s)) \nonumber \\ && \hspace{1.75in} - \; \sigma (X^{*}(s) , \nu_{\varepsilon } (s)) \big) - \sigma_x (X^{*}(s) , \nu^{*} (s) )) \, p (s) \, \big] \, d W(s) \nonumber \\
& & + \int_0^t S(t-s) \big[ \, \frac{1}{\varepsilon }\; \delta_{\varepsilon } \sigma (s) - \sigma_{\nu } ( X^{*}(s) , \nu^{*} (s) ) \, \nu (s) \, \big] \, ds
\nonumber \\
&& = \, \int_0^t S(t-s) \big[ \,  \tilde{b}_x (s , \varepsilon ) \, \eta_{\varepsilon} (s) +
( \, \tilde{b}_x (s , \varepsilon ) -  b_x ( X^{*}(s) , \nu^{*} (s) ) \, ) \, p (s) \, \big] \, ds \nonumber \\
& & + \int_0^t S(t-s) \big[ \, \frac{1}{\varepsilon }\; \delta_{\varepsilon } b (s) - b_{\nu } ( X^{*}(s) , \nu^{*} (s) ) \, \nu (s) \, \big] \, ds  \nonumber \\
& & +  \int_0^t S(t-s) \big[ \, \tilde{\sigma}_x (s , \varepsilon ) \, \eta_{\varepsilon} (s) +
( \, \tilde{\sigma}_x (s , \varepsilon ) -  \sigma_x ( X^{*}(s) , \nu^{*} (s) ) \, p (s) \, \big] \, d W(s)
\nonumber \\
& & + \int_0^t S(t-s) \big[ \, \frac{1}{\varepsilon }\; \delta_{\varepsilon } \sigma (s) - \sigma_{\nu } ( X^{*}(s) , \nu^{*} (s) ) \, \nu (s) \, \big] \, d W(s) ,
\end{eqnarray}
where \[ \delta_{\varepsilon } \sigma (s) = \sigma (X^{*}(s) , \nu_{\varepsilon } (s) ) - \sigma (X^{*}(s) , \nu^{*} (s) )  \] and
\[ \tilde{\sigma}_x (s , \varepsilon ) (y) = \int_0^1 \sigma_x (  X^{*}(s) + \theta ( X_{\varepsilon} (s) - X^{*} (s) ) , \nu_{\varepsilon } (s) ) (y) d\theta ,  \;\; y \in H . \]
Consequently, from (i) and  as in the proof of Lemma~\ref{lem:2nd estimate}, it follows that
\begin{eqnarray}\label{ineq1:lem3}
\mathbb{E}\, [ \, | \eta_{\varepsilon} (t) |^2 \, ] \leq C_6  \int_0^t  \mathbb{E}\,  [ \, | \eta_{\varepsilon} (s) |^2 \, ] ds + \rho ( \varepsilon  ),
\end{eqnarray}
for all $t \in [0, T] ,$ where
\begin{eqnarray}\label{eq:rho}
\rho ( \varepsilon  ) &=&  8 M T \,  \mathbb{E}\,  [ \, \int_0^T | \,  ( \, \tilde{b}_x (s , \varepsilon ) -  b_x ( X^{*}(s) , \nu^{*} (s) ) \, ) \, p (s) \, |^2  \; ds \, ] \nonumber \\
&&  +\; 8 M \, \mathbb{E}\,  [ \, \int_0^T || \, ( \, \tilde{\sigma}_x (s , \varepsilon ) -  \sigma_x ( X^{*}(s) , \nu^{*} (s) ) \, ) \, p (s) \, ||^2_{2}  \; ds \, ]
\nonumber \\ &&
+ \; 4 M T \, \mathbb{E}\,  [ \, \int_0^T | \, \frac{1}{\varepsilon }\; \delta_{\varepsilon } b (s) - b_{\nu } ( X^{*}(s) , \nu^{*} (s) ) \, \nu (s) \, |^2 \, ds \, ]
 \nonumber \\
&&  + \; 4 M \, \mathbb{E}\,  [ \, \int_0^T || \, \frac{1}{\varepsilon }\; \delta_{\varepsilon } \sigma (s) - \sigma_{\nu } ( X^{*}(s) , \nu^{*} (s) ) \, \nu (s) \, ||_2^2 \, ds \, ] .
\end{eqnarray}

But (i), (\ref{eq:est-p}) and the dominated convergence theorem give
\begin{eqnarray*}
&&  \hspace{-1cm} \mathbb{E}\,  [ \, \int_0^T | \,  ( \, \tilde{b}_x (s , \varepsilon ) -  b_x ( X^{*}(s) , \nu^{*} (s) ) \, ) \, p (s) |^2  \; ds \, ] \nonumber \\
&& \hspace{-0.5cm} = \, \mathbb{E}\,  [ \, \int_0^T   |\, \int_0^1 \big{(} \, b_x (  X^{*}(s) + \theta ( X_{\varepsilon} (s) - X^{*} (s) ) , \nu_{\varepsilon } (s) ) \nonumber \\
&& \hspace{2in} - b_x ( X^{*}(s) , \nu^{*} (s) ) \, \big{)} \, p (s) \, d\theta  |^2 \, ds \, ]
\nonumber \\
&& \hspace{-0.5cm} \leq
\int_0^T  \int_0^1  \mathbb{E}\,  [ \, |\,  \big{(} \, b_x (  X^{*}(s) + \theta ( X_{\varepsilon} (s) - X^{*} (s) ) , \nu_{\varepsilon } (s) )
\nonumber \\
&& \hspace{2in} - b_x ( X^{*}(s) , \nu^{*} (s) ) \, \big{)} \,  p (s) |^2 \, ] \, d\theta \, ds
\nonumber \\
&& \hspace{3.1in} \rightarrow 0 , \;\; \; \text{as}\; \; \varepsilon \rightarrow 0^{+} .
\end{eqnarray*}

Similarly we have
\begin{eqnarray}\label{eq:sigma-limit}
\mathbb{E}\,  [ \, \int_0^T || \, ( \, \tilde{\sigma}_x (s , \varepsilon ) -  \sigma_x ( X^{*}(s)) \, ) \, p (s) \, ||^2_{2}  \; ds \, ] \rightarrow 0 ,
\end{eqnarray}
as $\varepsilon \rightarrow 0^{+} .$

On the other hand, as done for (\ref{ineq4::difference of x}),
\begin{eqnarray}\label{ineq3:lem3}
&& \hspace{-2cm} \mathbb{E}\,  [ \, \int_0^T | \, \frac{1}{\varepsilon }\; \delta_{\varepsilon } b (s) - b_{\nu } ( X^{*}(s) , \nu^{*} (s) ) \, \nu (s) \, |^2 \, ds \, ]
\nonumber \\
&& \hspace{-1.5cm}
\leq \int_0^T \int_0^1  \mathbb{E}\,  \Big{[} \, \Big{|} \, \Big{(} \, b_{\nu } (  X^{*}(s) , \nu^{*} (s) + \theta ( \nu_{\varepsilon} (s) - \nu^{*} (s) ) ) \nonumber \\
&& \hspace{0.9in}
- b_{\nu } ( X^{*}(s) , \nu^{*} (s) ) \Big{)} \, \nu (s) \, \Big{|}^2 \, \Big{]} \, d\theta \, ds
\; \rightarrow 0 ,
\end{eqnarray}
if $\varepsilon \rightarrow 0^{+} ,$ by using (i) and the dominated convergence theorem.
Similarly,
\begin{eqnarray}\label{ineq4:lem3}
\mathbb{E}\,  [ \, \int_0^T || \, \frac{1}{\varepsilon }\; \delta_{\varepsilon } \sigma (s) - \sigma_{\nu } ( X^{*}(s) , \nu^{*} (s) ) \, \nu (s) \, ||_2^2 \, ds \, ] \rightarrow 0 ,
\end{eqnarray}
if $\varepsilon \rightarrow 0^{+} .$

Finally applying (\ref{eq:sigma-limit})--(\ref{ineq4:lem3}) in (\ref{eq:rho}) shows that
\[ \rho (\varepsilon ) \rightarrow 0 , \;\; \; \text{as}\; \; \varepsilon \rightarrow 0^{+} . \]
Hence from (\ref{ineq1:lem3}) and Gronwall's inequality we obtain
\[  \sup_{t \in [0 , T]} \mathbb{E} \, [ \, | \eta_{\varepsilon} (t) |^2 \, ]  \rightarrow 0 , \]
as $\varepsilon \rightarrow 0^{+} .$
\end{proof}

\bigskip

The following theorem contains our main variational equation, which is one of the main tools needed for deriving the maximum principle stated in Theorem~\ref{main thm}.
\begin{theorem}\label{thm:variational inequality}
We suppose that (i) and (ii) in Theorem~\ref{main thm} hold. For each $\varepsilon > 0 ,$ we have
\begin{eqnarray}\label{variational inequality}
&& \hspace{-2cm} J ( \nu_{ \varepsilon } ( \cdot ) ) - J ( \nu^{*} ( \cdot ) ) = \varepsilon \; \mathbb{E} \; [ \, \phi_x ( X^{*} (T) ) \, p(T) \, ]
\nonumber \\
&& + \; \varepsilon \; \mathbb{E} \; [ \, \int_0^T \ell_x (X^{*} (s) , \nu^{*} (s) ) \, p(s) \, ds \, ] \nonumber \\
&& + \; \mathbb{E} \; [ \, \int_0^T \big{(} \, \ell (X^{*} (s) , \nu_{\varepsilon } (s) ) - \ell (X^{*} (s) , \nu^{*} (s) ) \big{)} \, ds \, ]
+ o (\varepsilon ) .
\end{eqnarray}
\end{theorem}
\begin{proof}
We can write  $J ( \nu_{ \varepsilon } ( \cdot ) ) - J ( \nu^{*} ( \cdot ) )$ as
\begin{eqnarray}\label{eq1:var-ineq}
J ( \nu_{ \varepsilon } ( \cdot ) ) - J ( \nu^{*} ( \cdot ) ) = I_1 (\varepsilon ) + I_2 (\varepsilon ) ,
\end{eqnarray}
with
\begin{eqnarray*}
&& \hspace{-0.5cm} I_1 (\varepsilon ) = \mathbb{E} \; [ \, \phi ( X_{\varepsilon } (T) ) - \phi ( X^{*} (T) ) \, ]
\end{eqnarray*}
and
\begin{eqnarray*}
&& \hspace{-0.5cm} I_2 (\varepsilon ) = \mathbb{E} \; [ \, \int_0^T \big{(} \, \ell (X_{\varepsilon } (s) , \nu_{\varepsilon } (s) ) - \ell (X^{*} (s) , \nu^{*} (s) ) \big{)} \, ds \, ] .
\end{eqnarray*}

Note that with the help of our assumptions and by making use of Lemma~\ref{lem:3rd estimate}, Lemma~\ref{lem:2nd estimate}, Lemma~\ref{lem:1st estimate} and the dominated convergence theorem we deduce that
\begin{eqnarray*}
\hspace{-0.25cm} \frac{1}{\varepsilon } \, I_1 (\varepsilon ) &=& \frac{1}{\varepsilon } \, \mathbb{E} \, [ \, \int_0^1  \phi_x ( X^{*} (T) + \theta \, ( \, X_{\varepsilon } (T) - X^{*} (T) \, ) \, ( X_{\varepsilon } (T) - X^{*} (T) ) \, d\theta \, ] \nonumber \\
&& \hspace{-0.5cm}
= \mathbb{E} \, [ \, \int_0^1  \phi_x ( X^{* } (T) + \theta \, ( \, X_{\varepsilon } (T) - X^{*} (T) \, ) \, ( p (T) + \eta_{\varepsilon } (T) ) \, d\theta \, ] \nonumber \\
&& \hspace{-0.5cm} \rightarrow \mathbb{E} \, [ \, \phi_x ( X^{*} (T) ) \, p(T) \, ] , \;\; \; \text{as}\; \; \varepsilon \rightarrow 0^{+} .
\end{eqnarray*}
Hence
\begin{eqnarray}\label{eq4:var-ineq}
I_1 (\varepsilon ) = \varepsilon \; \mathbb{E} \; [ \, \phi_x ( X^{*} (T) ) \, p(T) \, ]
+ o (\varepsilon ) .
\end{eqnarray}
Similarly
\begin{eqnarray*}
&& \hspace{-0.5cm} \frac{1}{\varepsilon } \; I_2 (\varepsilon ) = \frac{1}{\varepsilon } \; \mathbb{E} \; [ \, \int_0^T \big{(} \, \ell (X_{\varepsilon } (s) , \nu_{\varepsilon } (s) ) - \ell (X^{*} (s) , \nu_{\varepsilon } (s) ) \big{)} \, ds \, ]
\nonumber \\
&& \hspace{1cm}
+ \; \frac{1}{\varepsilon } \; \mathbb{E} \; [ \, \int_0^T \big{(} \, \ell (X^{*} (s) , \nu_{\varepsilon } (s) ) - \ell (X^{*} (s) , \nu^{*} (s) ) \big{)} \, ds \, ]
\nonumber \\
&& = \mathbb{E} \; [ \, \int_0^T  \int_0^1 \ell_x (  X^{*}(s) + \theta ( X_{\varepsilon} (s) - X^{*} (s) ) , \nu_{\varepsilon } (s) ) \, ( p (s) + \eta_{\varepsilon } (s) ) \, d\theta \, ds \, ]
\nonumber \\
&& \hspace{2cm}
+ \; \frac{1}{\varepsilon } \; \mathbb{E} \; [ \, \int_0^T \big{(} \, \ell (X^{*} (s) , \nu_{\varepsilon } (s) ) - \ell (X^{*} (s) , \nu^{*} (s) ) \big{)} \, ds \, ] .
\end{eqnarray*}

On the other hand, applying Lemma~\ref{lem:3rd estimate}, Lemma~\ref{lem:2nd estimate}, Lemma~\ref{lem:1st estimate}, using the continuity and boundedness of $\ell_x$ in (i), (ii) and the dominated convergence theorem imply that
\begin{eqnarray*}
&& \hspace{-0.7cm} \mathbb{E} \; [ \, \int_0^T  \int_0^1 \ell_x (  X^{*}(s) + \theta ( X_{\varepsilon} (s) - X^{*} (s) ) , \nu^{* } (s) + \varepsilon \, \nu (s)) \, ( p (s) + \eta_{\varepsilon } (s) ) \, d\theta \, ds \, ]
\nonumber \\
&& \hspace{2.5in} \rightarrow
\mathbb{E} \; [ \, \int_0^T \ell_x (  X^{*}(s) , \nu^{*} (s) ) \, p (s)  ds \, ] .
\end{eqnarray*}
In particular we obtain
\begin{eqnarray}\label{eq6:var-ineq}
\hspace{-1.5cm} I_2 (\varepsilon ) &=& \varepsilon \;\mathbb{E} \; [ \, \int_0^T \ell_x (  X^{*}(s) , \nu^{*} (s) ) \, p (s)  ds \, ]
\nonumber \\
&& + \; \mathbb{E} \; [ \, \int_0^T \big{(} \, \ell (X^{*} (s) , \nu_{\varepsilon } (s) ) - \ell (X^{*} (s) , \nu^{*} (s) ) \big{)} \, ds \, ] + o (\varepsilon ) .
\end{eqnarray}

As a result the theorem follows from (\ref{eq1:var-ineq})--(\ref{eq6:var-ineq}).
\end{proof}

\bigskip

Let us next introduce an important variational inequality.
\begin{lemma}\label{lem:2ndvariational inequality}
Let hypotheses (i), (ii) in Theorem~\ref{main thm} hold. Let $(Y^{*} , Z^{*}) \equiv ( Y^{\nu^{*} (\cdot )},
Z^{\nu^{*} (\cdot )} )$ be the solution of BSEE~(\ref{adjoint-bse}) corresponding to the optimal pair $(X^{*}, \nu^{*}(\cdot)) .$ Then
\begin{eqnarray}\label{variational inequality2}
&& \hspace{-1.9cm} \varepsilon \; \mathbb{E} \; \big< \, Y^{*} (T) , p(T) \, \big>
+ \; \varepsilon \; \mathbb{E} \; [ \, \int_0^T \ell_x (X^{*} (s) , \nu^{*} (s) ) \, p(s) \, ds \, ] \nonumber \\
&& \hspace{-1.5cm} + \; \mathbb{E} \; \big[ \, \int_0^T \big( \,  \delta_{\varepsilon } \mathcal{H} (s) - \big<  \, \delta_{\varepsilon } b (s) \, , Y^{*}(s) \, \big> -  \big<  \, \delta_{\varepsilon } \sigma (s) \, , Z^{*}(s) \, \big>_2 \, \big) \, ds \, \big]
\geq o (\varepsilon ) ,
\end{eqnarray}
where
\begin{eqnarray*}
&& \delta_{\varepsilon } \mathcal{H} (s) =  \mathcal{H}
             ( X^{*} (s), \nu_{\varepsilon } (s), Y^{*} (s) , Z^{*} (s) )
 - \mathcal{H} ( X^{*} (s), \nu^{*} (s), Y^{*} (s) , Z^{*} (s) ) .
\end{eqnarray*}
\end{lemma}
\begin{proof}
Since $\nu^{*} ( \cdot )$ is an optimal control, then $J ( \nu_{ \varepsilon } ( \cdot ) ) - J ( \nu^{*} ( \cdot ) )  \geq 0 . $
Hence the result follows from (\ref{variational inequality}) and (\ref{def:Hamiltonian}).
\end{proof}

\bigskip

The following duality relation between (\ref{eq:p}) and (\ref{adjoint-bse}) is also needed in order to establish of proof of Theorem~\ref{main thm}.

\smallskip

\begin{lemma}\label{lem:formula of Y and p}
Under hypothesis (i) in Theorem~\ref{main thm}, we have
\begin{eqnarray}\label{Ito formula of Y and p}
&& \hspace{-2cm} \mathbb{E} \; \big< \, Y^{*} (T) , p(T) \, \big> = - \; \mathbb{E} \; [ \, \int_0^T \ell_x (X^{*} (s) , \nu^{*} (s) ) \, p(s) \, ds \, ] \nonumber \\
&& \hspace{2cm} + \; \mathbb{E} \; [ \, \int_0^T \big< \, b_{\nu } ( X^{*}(s) , \nu^{*} (s) ) \, \nu (s) \, , Y^{*}(s) \, \big> \, ds \, ]
\nonumber \\
&& \hspace{2cm} + \; \mathbb{E} \; [ \, \int_0^T \big< \, \sigma_{\nu } ( X^{*}(s) , \nu^{*} (s) ) \, \nu (s) \, , Z^{*}(s) \, \big>_2 \, ds \, ]  .
\end{eqnarray}
\end{lemma}
\begin{proof}
The proof is done by using Yosida approximation of the operator $A$ and It\^{o}'s formula for the resulting SDEs, and can be gleaned directly from the proof of Theorem~2.1 in \cite{[Tess96]}.
\end{proof}

\bigskip

We are now ready to establish (or complete in particular) the proof of Theorem~\ref{main thm}.

\bigskip

\noindent \begin{proof}[ {\bf Proof of Theorem~\ref{main thm}}]
Recall the BSEE~(\ref{adjoint-bse}):
\begin{eqnarray*}
\left\{ \begin{array}{ll}
             -\, d Y^{\nu (\cdot )} (t) = & \big(\, A^{*} \, Y^{\nu (\cdot )} (t) + \nabla_{x} \mathcal{H}
             ( X^{\nu (\cdot )}(t), \nu (t), Y^{\nu (\cdot )}(t) , Z^{\nu (\cdot )} (t) ) \,
             \big)\, dt \\& \hspace{1.70in} - Z^{\nu (\cdot )} (t) d W(t) , \; \; 0 \leq t < T,  \\
            \; \; \; \; Y^{\nu (\cdot )} (T) = & \nabla \phi (X^{\nu (\cdot )}(T)) .
         \end{array}
 \right.
\end{eqnarray*}
From Theorem~\ref{th:solution of adjointeqn} there exists a unique solution $( Y^{*},
Z^{*} )$ to it. Thereby it remains to prove (\ref{ineq1:main}).

Applying (\ref{variational inequality2}) and (\ref{Ito formula of Y and p}) gives
\begin{eqnarray}\label{eq1:proof of main theorem}
&& \hspace{-1.5cm}
\mathbb{E} \; \Big[ \, \int_0^T \Big( \,  \delta_{\varepsilon } \mathcal{H} (s) + \big<  \, \varepsilon \; b_{\nu } ( X^{*}(s) , \nu^{*} (s) ) \nu (s) - \delta_{\varepsilon } b (s) \, , Y^{*}(s) \, \big>
\nonumber \\ && \hspace{0.5cm}
+ \; \big<  \, \varepsilon \; \sigma_{\nu } ( X^{*}(s) , \nu^{*} (s) ) \nu (s) - \delta_{\varepsilon } \sigma (s) \, , Z^{*}(s) \, \big>_2 \, \Big) \, ds \, \Big]
\geq o (\varepsilon ).
\end{eqnarray}
But, as done for (\ref{ineq3:lem3}), by using the continuity and boundedness of $b_{\nu }$ in assumption (i) and the dominated convergence theorem, one can find that
\begin{eqnarray*}
&& \hspace{-0.75cm} \frac{1}{\varepsilon }\,  \mathbb{E} \; [ \, \int_0^T \big<  \, \varepsilon \; b_{\nu } ( X^{*}(s) , \nu^{*} (s) ) \nu (s) - \delta_{\varepsilon } b (s) \, , Y^{*}(s) \, \big> \, ds \\
&& = - \; \mathbb{E} \, \big{[} \, \int_0^T \big< Y^{*}(s) , \int_0^1  \Big{(} \, b_{\nu } (  X^{*}(s) , \nu^{*} (s) + \theta ( \nu_{\varepsilon} (s) - \nu^{*} (s) ) ) \nonumber \\ &&  \hspace{2in} - \, b_{\nu } ( X^{*}(s) , \nu^{*} (s) ) \Big{)} \, \nu (s) \,  \, d\theta  \,
 \big> \, ds \, \big{]} \; \rightarrow 0 ,
\end{eqnarray*}
as $\varepsilon \rightarrow 0^{+} .$
This means that
\begin{equation*}
\mathbb{E} \; [ \, \int_0^T \big<  \, \varepsilon \; b_{\nu } ( X^{*}(s) , \nu^{*} (s) ) \nu (s) - \delta_{\varepsilon } b (s) \, , Y^{*}(s) \, \big> \, ds \, ] = o (\varepsilon ) .
\end{equation*}
Similarly,
\begin{equation*}
\mathbb{E} \; [ \, \int_0^T \big<  \, \varepsilon \; \sigma_{\nu } ( X^{*}(s) , \nu^{*} (s) ) \nu (s) - \delta_{\varepsilon } \sigma (s) \, , Z^{*}(s) \, \big>_2 \, ds \, ] = o (\varepsilon ) .
\end{equation*}

Now by applying these two former identities in (\ref{eq1:proof of main theorem}) we deduce that
\begin{equation}\label{ineq2:main}
\mathbb{E} \; [ \, \int_0^T \delta_{\varepsilon } \mathcal{H} (s)  \, ds \, ]
\geq o (\varepsilon ) .
\end{equation}
Therefore, by dividing (\ref{ineq2:main}) by $\varepsilon$ and letting $\varepsilon \rightarrow 0^{+} ,$ the following inequality holds:
\begin{eqnarray*}
&&  \hspace{-1cm} \mathbb{E} \; [ \, \int_0^T \big{<} \nabla_{\nu }\mathcal{H} ( t , X^{*} (t) , \nu^{*} (t) ,
 Y^{*} (t) , Z^{*} (t) ) , \nu (t) \big{>}_{\mathcal{O}} \, dt \, ] \geq 0 .
\end{eqnarray*}

Finally, (\ref{ineq1:main}) follows by arguing, if necessary, as in \cite[P. 280]{[Be-book]} for instance.
\end{proof}

\bigskip

\bigskip

{\bf Acknowledgement.}
The author would like to thank the associate editor and anonymous referee(s) for their remarks and also for pointing out the recent work of Fuhrman et al., \cite{Fuh-Tes012}.

\fussy

\end{document}